\documentclass[a4paper, 12pt]{amsart}
\usepackage{amssymb,latexsym,amsmath,amsfonts,hhline}
\usepackage{pdfsync}
\usepackage{xcolor}
\input xy
\xyoption{all}

\def\grphp#1{$\xymatrix@R=10pt@C=10pt@M=0pt@L=2pt{#1}$}

\makeatletter
\renewcommand{\subsection}{\@startsection{subsection}{2}{0mm}{-2mm}{-2mm}{\bf\normalsize}}

\def\sbsnt#1{\subsection{#1}}
\makeatother

\newtheorem{formula}{}[section]
\newtheorem{definition}[formula]{Definition}
\newtheorem{corollary}[formula]{Corollary}
\newtheorem{remark}[formula]{Remark}
\newtheorem{lemma}[formula]{Lemma}
\newtheorem{theorem}[formula]{Theorem}

\def\thrm{\begin{theorem}}
\def\thrml#1{\begin{theorem}\label{#1}}
\def\ethrm{\end{theorem}}
\def\rmrk{\begin{remark}}
\def\rmrkl#1{\begin{remark}\label{#1}}
\def\ermrk{\end{remark}}
\def\dfntn{\begin{definition}}
\def\dfntnl#1{\begin{definition}\label{#1}}
\def\edfntn{\end{definition}}
\def\nmrt{\begin{enumerate}}
\def\enmrt{\end{enumerate}}
\def\tm#1{\item[{\rm (#1)}]}
\def\qtn{\begin{equation}}
\def\qtnl#1{\begin{equation}\label{#1}}
\def\eqtn{\end{equation}}
\def\lmm{\begin{lemma}}
\def\lmml#1{\begin{lemma}\label{#1}}
\def\elmm{\end{lemma}}
\def\crllr{\begin{corollary}}
\def\crllrl#1{\begin{corollary}\label{#1}}
\def\ecrllr{\end{corollary}}
\def\css{\begin{cases}}
\def\ecss{\end{cases}}
					
\def\proof{\noindent{\bf Proof}.\ }

\def\cE{{\mathcal E}}

\def\cK{{\mathcal K}}

\def\cX{{\mathcal X}}
\def\cY{{\mathcal Y}}
\def\cZ{{\mathcal Z}}

\def\fS{{\mathfrak S}}

\def\mZ{C}

\DeclareMathOperator{\aut}{Aut}

\DeclareMathOperator{\im}{im}

\DeclareMathOperator{\orb}{Orb}

\DeclareMathOperator{\poly}{poly}

\DeclareMathOperator{\reg}{Reg}

\DeclareMathOperator{\sym}{Sym}
\DeclareMathOperator{\WL}{WL}

\def\eprf{\hfill$\square$}

\def\qaq{\quad\text{and}\quad}
\def\qoq{\quad\text{or}\quad}

\newcommand{\grp}[1]{\langle {#1}\rangle}

\title[Recognising and testing isomorphism of Cayley graphs]{Recognizing and testing isomorphism of Cayley graphs over an abelian group of order $4p$  in polynomial time}

\author[R. Nedela]{Roman Nedela}
\address[R. Nedela]{Faculty of Applied Sciences, University of West Bohemia, Technick\'a 8, Pilsen, Czech Republic}
\email[R. Nedela]{nedela@savbb.sk}
\author[I. Ponomarenko]{Ilia Ponomarenko}
\address[I. Ponomarenko]{V. A. Steklov Institue of Mathematics, Russian Academy of Sciences, Sankt Petersburg, Russia}
\email[I. Ponomarenko]{inp@pdmi.ras.ru}

\begin{document}

\begin{abstract}
We construct a polynomial-time algorithm that given a graph $X$ with $4p$ vertices ($p$ is prime), finds (if any) a  Cayley representation of $X$ over the group $C_2\times C_2\times C_p$. This result, 
together with the known similar result for circulant graphs,
 shows that recognising and testing isomorphism of Cayley graphs over an abelian group of order $4p$ can be done in polynomial time.
\end{abstract}

\maketitle

\section{Introduction}

Under a {\it Cayley representation} of a graph $X$ over a group~$G$, we mean a graph isomorphism from $X$ to a Cayley graph over~$G$. Two such  representations  are {\it equivalent} if the images are {\it Cayley isomorphic}, i.e., there is a group automorphism of $G$ which is a graph isomorphism between the images. In the present paper, we consider a special case of the following computational problem (below all the groups and graphs are assumed to be finite). \medskip
 
\noindent {\bf Problem CRG.} {\it Given a group $G$ and a graph $X$, find a full set of non-equivalent Cayley representations of $X$ over $G$.}\medskip
 
Here we assume that the group $G$ is given explicitly by the multiplication table, and the graph $X$ is given by a binary relation. The output of an algorithm solving the problem is represented by a set of bijections~$f$ from the vertex set of $X$ onto $G$ such that $(G_{right})^{f^{-1}}$ is a regular subgroup of the group $\aut(X)$, where $G_{right}\le\sym(G)$ is the group induced by right multiplication of $G$. Using the Babai argument in \cite{B77}, one can establish a one-to-one correspondence between the regular subgroups of $\aut(X)$ and Cayley representations of $X$. \medskip

In general, the Problem CRG seems to be very hard. Even the question whether the output is empty or not, leads to the Cayley recognition problem asking whether a given graph is isomorphic to a Cayley graph over the group $G$. Not too much is known about the computational complexity of this problem. There are two other related problems.\medskip

\noindent {\bf Problem CGCI.} {\it Given a group $G$, test whether two Cayley graphs over $G$ are Cayley isomorphic.}\medskip

\noindent {\bf Problem CGI.} {\it Given a group $G$, a Cayley graph over $G$, and an arbitrary graph, test whether these two graphs are isomorphic.}\medskip

Note that the subproblem of the Problem CGI, in which both input graphs are Cayley graphs over $G$, is equivalent to the Problem CGCI whenever $G$ is a CI-group (this fact can be considered as the definition of the CI-group.)\medskip

Suppose we are restricted to a family of Cayley graphs for which the Problem CGCI can be solved efficiently; for instance, one can take~$G$ to be a group generated by a set of at most constant size. Then one can see that the Problem CGI is polynomial-time reducible to the Problem CRG. In general case, the reduction can be done in polynomial time in the order of the group~$\aut(G)$.\medskip

In  \cite{EP04}, a polynomial-time algorithm for the Problem CRG was constructed for the case where $G$ is a cyclic group. Up to now, this is the only published result solving the Problem CRG for an infinite class of groups. It is quite natural to look for an extension of that result to other classes of abelian groups. The main result of the present paper (Theorem~\ref{100916a}) does it for abelian groups of order $4p$, where $p$ is a prime. In view of the above discussion, the Problems CRG and CGI are equivalent in this case.\medskip

Now we are ready to present the main results of the paper.

\thrml{100916a}
For an abelian group $G$ of order $n=4p$ with prime $p$, the Problems CRG and CGI can be solved  in time $\poly(n)$.
\ethrm

There are exactly two non-isomorphic abelian groups of order $4p$: the cyclic group $C_{4p}$ and the group $E_4\times C_p$, where $E_4=C_2\times C_2$ is the Klein group. In the former case, Theorem~\ref{100916a} follows from \cite{EP04}. In the latter case, $G$ is a CI-group \cite[Theorem~1.2]{KM09} and hence every graph has at most one Cayley representation over $G$ (up to equivalence). Thus Theorem~\ref{100916a} is an immediate consequence of the following theorem.

\thrml{030916a}
Given a graph $X$ with $n=4p$ vertices ($p$ is a prime), one can test in time $\poly(n)$ whether $X$ is isomorphic to a Cayley graph over the group $G=E_4\times\mZ_p$ and (if so) find a Cayley representation of~$X$ over $G$ within the same time.
\ethrm

Let us outline the proof of Theorem~\ref{030916a}. At the first step, we use  the Weisfeiler-Leman algorithm to construct the coherent configuration~$\cX$ associated with the graph~$X$ (for the exact definitions, see Section~\ref{210217a}). Then
$$
K:=\aut(\cX)=\aut(X).
$$
Therefore $X$ is a Cayley graph over~$G$ if and only if $\cX$ is a Cayley scheme over~$G$. This reduces our problem to finding the set 
$\reg(K,G)$ consisting of all semiregular groups $H\leq K$ isomorphic to $G$, where the group $K$ is not ``in hand''.\medskip 

At this point, we use the classification of Schur rings (and hence Cayley schemes) over the group~$G$ obtained in~\cite{EKP13}\footnote{Independently, such a classification have recently be obtained in~\cite{L16}.}. This 
enables us to find in time $\poly(n)$ a larger coherent configuration $\cX'$ 
such that
\nmrt
\tm{a} $\reg(K,C_p)=\reg(K',C_p)$, where $K'=\aut(\cX')$,
\tm{b} $|K'|=\poly(n)$ or $K'\cong\sym(p')^m$, where $p'm=n$ and $m\le 4$.
\enmrt
Now, if the group $K'$ has small order, then in view of statement~(a), one can easily find the set $\reg(K,C_p)$ by brute force; this is a part of the Main Subroutine described in Section~\ref{14102016}. In the remaining case, statement~(b) implies that $\cX'$ is a coherent configuration of a special type studied in Section~\ref{14102016t}. This fact is used in the Main Subroutine for computing the set $\reg(K,C_p)$. At the final step, we only need to test whether there exists a group belonging to $\reg(K,C_p)$, which can be extended to a regular subgroup of $K$. This is done in Section~\ref{100916g}.

\section{Coherent configurations}\label{210217a}

In this section, we collect some notation and then compile basic definitions and facts concerning coherent configurations. In our presentation, we follow~\cite{EP09}.

\sbsnt{Notation.} Throughout the paper, $\Omega$ denotes a finite set
of cardinality $n\geq 1$. The diagonal of the Cartesian product $\Omega\times\Omega$ is denoted by~$1_\Omega$. For a set $T\subseteq 2^\Omega$, we denote by $T^\cup$ the set of all unions of the elements of~$T$.  For a set $S\subseteq 2^{\Omega\times\Omega}$,  we set $S^*=\{s^*:\ s\in S\}$, where $s^*=\{(\beta,\alpha):\ (\alpha,\beta)\in s\}$. For a point $\alpha\in\Omega$, we set $\alpha S=\cup_{s\in S}\alpha s$, where $\alpha s=\{\beta\in\Omega:\ (\alpha,\beta)\in s\}$. For $r,s\subseteq\Omega\times\Omega$ we set $r\cdot s=\{(\alpha,\gamma):\ (\alpha,\beta)\in r,\ (\beta,\gamma)\in s$ for some $\beta\in\Omega\}$.\medskip

For an equivalence relation $E$ on $\Omega$, we denote by $\Omega/E$ the set of the classes of~$E$. If, in addition, $r\subseteq \Omega\times\Omega$, then we set
$$
r_{\Omega/E}=\{(\Lambda,\Delta)\in\Omega/E\times\Omega/E: r_{\Lambda,\Delta}\ne\varnothing\}
$$
where $r_{\Lambda,\Delta}=r\cap (\Lambda\times\Delta)$. We also 
put $r_\Lambda=r_{\Lambda,\Lambda}$. \medskip

The group of all permutations of $\Omega$ is denoted by $\sym(\Omega)$. Given a group $K\le\sym(\Omega)$ and  a $K$-invariant set $\Delta\subseteq\Omega$, the restrictions $k^\Delta$ of $k\in K$ to $\Delta$ form a subgroup of $\sym(\Delta)$ denoted by $K^\Delta$.\medskip

A bijection $f:\Omega \to \Omega'$, $\alpha\mapsto \alpha^f$, naturally defines a bijection $r\mapsto r^f$ from the relations on~$\Omega$ onto the relations on~$\Omega'$ and a group isomorphism $g\mapsto g^f$ from $\sym(\Omega)$ onto $\sym(\Omega')$. For an equivalence relation $E$ on~$\Omega$ the bijection $f$ induces a bijection $f_{\Omega/E}:\Omega/E\to \Omega'/E'$ where $E'=E^f$.

\sbsnt{Main definitions.}
Let $S$ be a partition of the set $\Omega\times\Omega$. A pair $\cX=(\Omega,S)$ is called a {\it coherent configuration} on $\Omega$ if $1_\Omega\in S^\cup$, $S^*=S$, and given $r,s,t\in S$, the number
$$
c_{rs}^t=|\alpha r\cap \beta s^*|
$$
does not depend on the choice of $(\alpha,\beta)\in t$. The elements of $\Omega$, $S$, $S^\cup$, and the numbers $c_{rs}^t$ are called the {\it points}, the {\it basis relations}, the {\it relations}, and the {\it intersection numbers} of~$\cX$, respectively. The numbers $|\Omega|$ and $|S|$ are called the {\it degree} and the {\it rank} of~$\cX$.\medskip

The coherent configuration $\cX$ is said to be {\it trivial} if the set $S$  consists of the reflexive relation $1_\Omega$ and (if $n>1$) the complement of it in $\Omega\times\Omega$, and is called {\it complete} if
every element of $S$ is a singleton. Observe that the rank of $\cX$ is at most two in the former case, and equals $n^2$ in the latter case.
\medskip

The set of all equivalence relations $E\in S^\cup$ is denoted by~$\cE=\cE(\cX)$. The coherent configuration is said to be {\it primitive} if the only elements of~$\cE$ are the trivial equivalence relations $1_\Omega$ and $\Omega\times\Omega$.

\sbsnt{Fibers.}
Denote by $\Phi=\Phi(\cX)$ the set of all $\Delta\subseteq\Omega$ such that $1_\Delta\in S$. Then the set $\Omega$ is the disjoint union of the elements of $\Phi$ called the {\it fibers} of~$\cX$. Moreover, for each $r\in S$ there exist uniquely determined fibers $\Delta$ and $\Lambda$ such that $r\subseteq\Delta\times\Lambda$. Thus the set $S$ is the disjoint union of the
sets
$$
S_{\Delta,\Lambda}=\{s\in S:\ s\subseteq\Delta\times\Lambda\}.
$$
Note that $1_\Delta\in rr^*$ and hence the number 
$$
n_r:=|\alpha r|=c_{rr^*}^{1_\Delta}
$$
does not depend on $\alpha\in\Delta$. It is called the {\it valency} of~$r$. For any $T\in S^\cup$, we set $n_T$ to be the sum of all valences $n_t$, where $t$ runs over the basis relations of~$\cX$  that are contained in~$T$.\medskip 

A coherent configuration $\cX$ is said to be {\it homogeneous}  if $1_\Omega\in S$. In this case, $n_r=n_{r^*}$ for all $r\in S$. Observe that a primitive coherent configuration is always homogeneous. A homogeneous coherent configuration which is not primitive is said to be {\it imprimitive}.
One can see that every {\it commutative} coherent configuration, i.e., one with $c_{rs}^t=c_{sr}^t$ for all $r,s,t$, is homogeneous.

\sbsnt{Restrictions and quotients.}
Let $E\in\cE$ be an equivalence relation. For any its class $\Delta$, denote by $S_\Delta$ the set of all nonempty relations~$s_\Delta=s\cap\Delta^2$ with $s\in S$. Then the pair
$$
\cX_\Delta=(\Delta,S_\Delta)
$$ 
is a coherent configuration called the {\it restriction} of~$\cX$ to the set~$\Delta$. In the special case where $E$ is the union of $\Lambda\times\Lambda$, $\Lambda\in\Phi$, and $\Delta\in\Phi$, the restriction $\cX_\Delta$ is called the {\it homogeneous component} of~$\cX$.\medskip
 
Let $\cX$ be a homogeneous coherent configuration. Denote by $S_{\Omega/E}$ the set of all nonempty relations $s_{\Omega/E}$, $s\in S$. Then the pair
$$
\cX_{\Omega/E}=(\Omega/E,S_{\Omega/E})
$$ 
is a coherent configuration called the {\it quotient} of $\cX$ modulo~$E$. Let $F\in\cE$ be an equivalence relation contained in~$E$. For a class $\Delta\in\Omega/E$, the equivalence relation $F_\Delta=F\cap \Delta^2$ on~$\Delta$ obviously belongs to the set $\cE(\cX_\Delta)$. On the other hand, the set $\Delta/F_\Delta$ is a class of the equivalence relation $E_{\Omega/F}$  being the union of basis relation of the quotient $\cX_{\Omega/F}$. We have the following commuting rule:
\qtnl{15102016}
(\cX_\Delta)_{\Delta/F_\Delta}=(\cX_{\Omega/F})_{\Delta/{F_\Delta}}.
\eqtn
The coherent configuration defined in \eqref{15102016} is denoted $\cX_{\Delta/F}$.

\sbsnt{Isomorphisms.}
Two coherent configurations are called {\it isomorphic} if there exists a bijection between their point sets that induces a bijection between their sets of basis relations. Each such bijection is called an {\it isomorphism} between these two configurations. The group of all isomorphisms of a coherent configuration $\cX$ to itself contains a normal subgroup
$$
\aut(\cX)=\{f\in\sym(\Omega):\ s^f=s,\ s\in S\}
$$
called the {\it automorphism group} of~$\cX$, where $s^f$ is the set of all pairs $(\alpha^f,\beta^f)$ with $(\alpha,\beta)\in s$. Thus by definition, $\aut(\cX)$ is the intersection of the automorphism groups  of the basis relations of $\cX$. It is easily seen that if $\Delta\in\Phi^\cup$ and $E\in\cE$, then
$$
\aut(\cX)^\Delta\le\aut(\cX_\Delta)\qaq \aut(\cX)^{\Omega/E}\le\aut(\cX_{\Omega/E}),
$$
where $\aut(\cX)^\Delta$ and $\aut(\cX)^{\Omega/E}$ are the permutations groups induced by the actions of the setwise stabilizer of~$\Delta$ in~$\aut(\cX)$ and of the group $\aut(\cX)$ on the sets $\Delta$ and $\Omega/E$, respectively.

\sbsnt{Direct sum.}
Let $\cX_1=(\Omega_1,S_1)$ and $\cX_2=(\Omega_2,S_2)$ be coherent configurations. Denote by  $\Omega_1\sqcup\Omega_2$ the disjoint union of the sets~$\Omega_1$ and~$\Omega_2$. Further, denote by $S_1\boxplus S_2$ the disjoint union of the set $S_1\sqcup S_2$ and the set of all Cartesian products $\Delta\times\Lambda$, where $\Delta\in \Phi(\cX_i)$ and
$\Lambda\in\Phi(\cX_{3-i})$, $i=1,2$. Then the pair
$$
\cX_1\boxplus\cX_2=(\Omega_1\sqcup\Omega_2,S_1\boxplus S_2)
$$
is a coherent configuration called the {\it direct sum} of~$\cX_1$ and~$\cX_2$. It is easily seen that
\qtnl{290916a}
\aut(\cX_1\boxplus\cX_2)=\aut(\cX_1)\times\aut(\cX_2),
\eqtn
where the direct product on the rihgt-hand side acts on the set $\Omega_1\sqcup\Omega_2$.

\sbsnt{Tensor product.}
Let $\cX_1=(\Omega_1,S_1)$ and $\cX_2=(\Omega_2,S_2)$ be coherent configurations.
Set $S_1\otimes S_2=\{s_1\otimes s_2:\ s_1\in S_1,\ s_2\in S_2\}$, where $s_1\otimes s_2$ is
the relation on $\Omega_1\times\Omega_2$ consisting of all 
pairs $((\alpha_1,\alpha_2),(\beta_1,\beta_2))$
with $(\alpha_1,\beta_1)\in s_1$ and $(\alpha_2,\beta_2)\in s_2$. Then the pair
$$
\cX_1\otimes\cX_2=(\Omega_1\times\Omega_2,S_1\otimes S_2)
$$
is a coherent configuration called the {\it tensor product} of~$\cX_1$ and~$\cX_2$. It is easily seen that
$$
\aut(\cX_1\otimes\cX_2)=\aut(\cX_1)\times\aut(\cX_2),
$$
where the direct product on the rihgt-hand side acts on the set $\Omega_1\times\Omega_2$.\medskip

Let $\cX$ be a commutative coherent configuration. We say that an equivalence relation $E\in\cE$ has a {\it complement} with respect to tensor product if there exists an equivalence relation $F\in\cE$ such that 
\qtnl{111216a}
E\cap F=1_\Omega\qaq E\cdot F=\Omega\times\Omega,
\eqtn
and $\cX$ is isomorphic to $\cX_{\Omega/E}\otimes\cX_{\Omega/F}$. It should be noted that for a fixed equivalence relations $E$ and $F$ satisfying conditions~\eqref{111216a}, the isomorphism $\cX\to \cX_{\Omega/E}\otimes\cX_{\Omega/F}$ exists if and only if $|\Omega/E|\cdot|\Omega/F|=n^2$, see \cite[Theorem 2.2]{EP04}.

\sbsnt{Generalized wreath product.}
Let $\cX$ be a homogeneous coherent configuration, and let $E$ and $F$ be equivalence relations belonging to the set $\cE$. We say that $\cX$ is the {\it $F/E$-wreath product} if $E\subseteq F$ and for each $r\in S$,
$$
r\cap F=\varnothing\qquad\Rightarrow\qquad r=\bigcup_{(\Delta,\Lambda)\in r_{\Omega/E}}\Delta\times\Lambda.
$$
When the equivalence relations are not relevant, we also say that $\cX$ is a {\it generalized wreath product}. It is said to be {\it trivial} if $E=1_\Omega$ or $F=\Omega\times\Omega$. The standard wreath product is obtained as a special case of the generalized wreath product with~$F=E$ (see also~\cite[p.45]{W76}). Finally, we say that an equivalence relation $E\in\cE$ has a {\it complement} with respect to (generalized) wreath product if there exists an equivalence relation $F\in\cE$ such that $\cX$ is the $F/E$-wreath product. 

\sbsnt{Algorithms.} 
From the algorithmic point of view, a coherent configuration $\cX$ on $n$ points is given by the set $S$ of its basis relations. In this representation, one can check in time $\poly(n)$ whether $\cX$ is homogeneous, commutative, etc. Moreover, within the same time one can list the fibers of $\cX$ and construct the restriction $\cX_\Delta$ for any $\Delta\in\Phi^\cup$, and can also find a nontrivial equivalence relation $E\in S^\cup$ (if it exists) and construct the quotient of~$\cX$ modulo~$E$.\medskip

\section{Quasitrivial coherent configurations.}\label{14102016t} 
Let $\cX=(\Omega,S)$ be a coherent configuration. It is said to be {\it quasitrivial} if 
$$
\aut(\cX)^\Delta=\sym(\Delta)\quad\text{for all}\ \Delta\in\Phi,
$$ 
where $\Phi=\Phi(\cX)$. In particular, the restriction  $\cX_\Lambda$ with $\Lambda\in\Phi^\cup$ of a quasitrivial configuration is quasitrivial as well. Thus, every homogeneous component of $\cX$ is trivial and homogeneous quasitrivial coherent configurations are exactly trivial ones.

\lmml{090916a}
The coherent configuration $\cX$ is quasitrivial if and only if for any two its fibers $\Delta$ and $\Lambda$, the set $S_{\Delta,\Lambda}$ is a singleton or contains exactly two elements, one of which is a bijection $f:\Delta\to\Lambda$.\footnote{This bijection is treated as the binary relation coinciding with the graph of~$f$.}
\elmm 
\proof To prove the ``if'' part, we make use of \cite[Lemma~9.4]{E03} implying that if $S_{\Delta,\Lambda}$ contains a bijection $f:\Delta\to\Lambda$, then the restriction map 
$$
\aut(\cX)\to\aut(\cX_{\Omega\setminus\Lambda})
$$ 
is an isomorphism. Using the induction on $|\Phi|$, one can reduce the general case to the case where $S_{\Delta,\Lambda}$ is a singleton for all $\Lambda\in\Phi$ other than~$\Delta$. But then it is easily seen that $\cX=\cX_\Delta\boxplus\cX_{\Omega\setminus\Delta}$. By formula~\eqref{290916a}, this implies that
$$
\aut(\cX)^\Delta=\aut(\cX_\Delta)=\sym(\Delta),
$$
where the last equality is true, because the coherent configuration $\cX_\Delta$ is of rank $|S_{\Delta,\Delta}|\le 2$.\medskip

To prove the ``only if'' part, we assume without loss of generality that $\Delta$ and $\Lambda$ are distinct non-singletons. Then the group $K:=\aut(\cX)^{\Delta\cup\Lambda}$ is the subdirect product of the groups $\sym(\Delta)$ and $\sym(\Lambda)$, i.e.,
$$
K=\{(g,h)\in \sym(\Delta)\times\sym(\Lambda):\ \varphi(g)=\psi(h)\}.
$$
where $\varphi:\sym(\Delta)\to M$ and $\psi:\sym(\Lambda)\to M$ are  epimorphisms to a suitable group $M$. Each of the groups $\ker(\varphi)$ and $\ker(\psi)$ is a normal subgroup of a $2$-transitive group and hence is transitive or trivial. If one of them is nontrivial, then  $S_{\Delta,\Lambda}$ is obviously a singleton and we are done. Now let $\ker(\varphi)=\ker(\psi)= 1$, i.e., $\varphi$ and $\psi$ are group isomorphisms. Then $M\cong K$ and hence 
$$
\sym(\Delta)\cong K\cong\sym(\Lambda).
$$ 
Consequently, $|\Delta|=|\Lambda|$. Denote the latter number by $d$. Note that by the assumption, $d\ge 2$. Therefore,
$$
K_\alpha\cong\sym(d-1)\quad\text{for all}\ \alpha\in\Delta\cup\Lambda
$$
and hence every non-singleton orbit of the group $K_\alpha\le\sym(\Delta\cup\Lambda)$ is of cardinality at least $d-1$. Therefore, there exists a relation $f\in S_{\Delta,\Lambda}$ such that $n_f=1$.
In other words, $f:\Delta\to\Lambda$ is a bijection taking $\delta\in\Delta$ to a unique point of the set $\delta f$. Finally in view of \cite[Corollary~13, p.~86]{W76},  we have
$$
|S_{\Delta,\Lambda}|\le \frac{1}{2}(|S_{\Delta,\Delta}|+|S_{\Lambda,\Lambda}|)=2
$$
which completes the proof.\eprf\medskip

Let the  coherent configuration $\cX$ be quasitrivial. We define a binary relation $\sim$ on the set $\Phi$ by setting 
\qtnl{111216c}
\Delta\sim\Lambda\quad\Leftrightarrow\quad |S_{\Delta,\Lambda}|=2.
\eqtn
This relation is obviously reflexive and symmetric. Moreover, assume that $|S_{\Lambda,\Gamma}|=2$ for some $\Gamma\in\Phi$. Then by Lemma~\ref{090916a}, the set $S_{\Lambda,\Gamma}$ contains a  bijection $g:\Lambda\to\Gamma$. It is easily seen that the superposition $fg:\Delta\to\Gamma$ coincides with $f\cdot g$ and hence belongs to $S_{\Delta,\Gamma}$. It follows that $\Delta\sim\Gamma$. Thus, $\sim$ is an equivalence relation.\medskip

Denote by $\Phi_1,\ldots,\Phi_m$ the classes of the equivalence relation $\sim$. For each $i=1,\ldots,m$, set
$$
\Delta_i=\Delta_{i1}\cup\cdots\cup\Delta_{im_i}
$$
where $\Delta_{i1},\ldots,\Delta_{im_i}$ are  the fibers belonging to the class~$\Phi_i$. Then formula~\eqref{111216c} immediately implies that 
\qtnl{091016a}
\cX={\overset{m}{\underset{i=1}{\boxplus}}}\cX_{\Delta_i}.
\eqtn
By formula~\eqref{290916a}, this enables us to find the automorphism group of a quasitrivial coherent configuration. Let $f_{ij}:\Delta_{i1}\to\Delta_{ij}$ be a bijection belonging to the set $S_{\Delta_{i1},\Delta_{ij}}$ (note that $f_{ij}$ is uniquely determined unless $|\Delta_{i1}|=2$).

\thrml{090916b}
Let $\cX$ be a quasitrivial coherent configuration. Then in the above notation, 
$$
\aut(\cX)\cong\prod_{i=1}^m\sym(\Delta_{i1}),
$$
\noindent Moreover, $g\in\aut(\cX)$ if and only if for every $i=1,\ldots,m$ there exists $g_i\in\sym(\Delta_{i1})$ such that
\qtnl{101016a}
g^{\Delta_{ij}}=\css
g_i,                     &\text{if $j=1$,}\\
f_{ij}^{-1}g_if_{ij}^{}  &\text{if $j>1$}.\\
\ecss
\eqtn
\ethrm 
\proof The second statement is an immediate consequence of the first one, which follows from formulas~\eqref{290916a} and ~\eqref{091016a}.\eprf

\section{Cayley schemes}
\sbsnt{General facts.} Let $G$ be a group. A coherent configuration~$\cX$ is called {\it a Cayley scheme} over~$G$ if 
$$
\Omega=G\qaq G_{right}\le\aut(\cX).
$$ 
In this case, the coherent configuration $\cX$ is homogeneous, and commutative if the group $G$ is abelian. We note that each basis relation of~$\cX$ is the arc set of a Cayley graph on~$G$. If the group $G$ is cyclic, then $\cX$ is said to be a {\it circulant} scheme.\medskip

For every equivalence relation $E\in\cE$, denote by $H=H_E$ the class of $E$ containing the identity element of~$G$. Then $H$ is a subgroup of~$G$ (isomorphic to the setwise stabilizer of $H$ in $G_{right}$) and the classes of~$E$ equal the right $H$-cosets of~$G$. In particular, 
$$
n_E=|H| \qaq|\Omega/E|=|G/H|.
$$ 
In what follows, we set $\cX_{G/H}=\cX_{\Omega/E}$. In the case where $\cX$ is the $F/E$-wreath product for some $E,F\in\cE$, we also say that $\cX$ is the $U/L$-wreath product, where $U=H_F$ and $L=H_E$.\medskip

We set
$$
\fS=\fS(\cX)=\{H_E\le G:\ E\in\cE(\cX)\}.
$$
For a group $H\in\fS$, we denote by $E_H$ the equivalence relation $F$ on~$G$, for which $H=H_F$. Thus, $E_{H_E}=E$ and $H_{E_H}=H$.  The following statement is a consequence of \cite[Theorem~5.6]{EP12}.

\thrml{121016O} 
Let $\cX$ be a Cayley scheme over an abelian group~$G$. Suppose that $\cX$ is the $U/L$-wreath product for some $U,L\in\fS$ such that 
\qtnl{121216a}
\aut(\cX_{G/L})^{U/L}=\aut(\cX_{U})^{U/L}.
\eqtn
Then $\aut(\cX)^{U}=\aut(\cX_{U})$ and $\aut(\cX)^{G/L}=\aut(\cX_{G/L})$.
\ethrm
\proof Set $\Delta_U=\aut(\cX_{U})$ and $\Delta_0=\aut(\cX_{G/L})$. Then according to \cite[Subsection~5.2]{EP12}, equality~\eqref{121216a} enables us to define the canonical generalized wreath product 
$$
K=\Delta_U\wr_{U/L}\Delta_0. 
$$
Since $\aut(\cX)^U\le\Delta_U$ and $\aut(\cX)^{G/L}\le\Delta_0$, we conclude by \cite[Corollary~5.5]{EP12} that $\aut(\cX)\le K$. The reverse inclusion follows from \cite[Corollary~5.4]{EP12}. Thus, $\aut(\cX)=K$. This implies that
$$
\aut(\cX)^{U}=K^U=\Delta_U=\aut(\cX_{U})
$$
and similarly $\aut(\cX)^{G/L}=\aut(\cX_{G/L})$.\eprf\medskip

A Cayley scheme~$\cX=(G,S)$ is said to be {\it cyclotomic} if there exists a group $K\le\aut(G)$ such that
$$
S=\orb(G_{right}K,G\times G).
$$ 
The cyclotomic scheme $\cX$ is called {\it proper} (respectively, {\it trivial}) if $K$ is a proper subgroup (respectively, the identity subgroup) of the group~$\aut(G)$. Note that for every cyclotomic scheme $\cX$, the set $\fS(\cX)$ contains all characteristic subgroups of $G$. 

\sbsnt{Cayley schemes over $E_4\times\mZ_p$.} The following theorem was proved in~\cite[Subsection~6.2]{EKP13} in the language of S-rings. Though the statement of the theorem  was not formulated explicitly there,  a careful analyses of the proof reveal cases (1) and (2) below, see \cite{L16} as well. 

\thrml{030916c}
Let $\cX$ be a Cayley scheme over the group~$G=E_4\times\mZ_p$. Then one of the following statements hold:
\nmrt
\tm{1} $\cX$ is trivial or cyclotomic,
\tm{2} $\cX$ is a nontrivial tensor or generalized wreath product.
\enmrt
\ethrm


\crllrl{12102016d}
Let $\cX$ be a Cayley scheme over the group~$G$ and $H$ a minimal subgroup in the set~$\fS=\fS(\cX)$. Then the Cayley scheme $\cX_H$ is trivial or isomorphic to a proper cyclotomic scheme over~$\mZ_p$.
\ecrllr
\proof Clearly, $|H|$ is a divisor of~$4p$ other than~$1$. Next, by the minimality of $H$ the scheme $\cX_H$ is primitive. If $|H|=4p$, then 
$\cX_H=\cX$ and $\fS=\{1,H\}$. This implies that the scheme $\cX$ is neither cyclotomic nor a nontrivial tensor or generalized wreath product. By Theorem~\ref{030916c}, we conclude that $\cX$ and hence $\cX_H$ is trivial. Now let $|H|=4$. Then $\cX_H$ is a primitive scheme over $E_4$ and hence is trivial. In the other three cases, the group~$H$ is a cyclic group of order $2$, $p$, or $2p$. Thus, the required statement follows from the Schur-Wielandt theory \cite[Corollary~3.2, Theorem~3.4]{MP}.\eprf\medskip

In \cite[p.423]{KP}, the automorphism groups of S-rings over a cyclic group $\mZ_{pq}$ with primes $q\ne p$, were completely classified. This classification shows that for every  Cayley scheme $\cY$ over the group $M=C_{2p}$ and for every group $N\le\fS(\cY)$, we have
\qtnl{121016t}
\aut(\cY)^N=\aut(\cY_N)\qaq \aut(\cY)^{M/N}=\aut(\cY_{M/N}).
\eqtn
Note that relations in \eqref{121016t} obviously hold also if $M=E_4$, or  if $M=C_p$. 

\lmml{12102016e}
Let $\cX$ be a Cayley scheme over the group~$G=E_4\times\mZ_p$. Suppose that the groups $U,L\in\fS$ are such that either $\cX=\cX_U\otimes\cX_L$, or $\cX$ is the $U/L$-wreath product. Then given $H\in\fS$,
\qtnl{121016a1}
\aut(\cX)^H=\aut(\cX_H)\quad\text{for all}\ H\leq L,
\eqtn
and 
\qtnl{121016a2}
\aut(\cX)^{G/H}=\aut(\cX_{G/H})\quad\text{for all}\ H\geq U.
\eqtn
\elmm
\proof Let us prove formula~\eqref{121016a1};  formula~\eqref{121016a2} can be proved in a similar way. Assume first that $\cX=\cX_U\otimes\cX_L$. Then by the first part of formula~\eqref{121016t} for $\cY=\cX_L$, $M=L$, and $N=H$, we have
$$
\aut(\cX)^H=(\aut(\cX_U)\otimes\aut(\cX_L))^H=\aut(\cX_L)^H=\aut(\cX_H),
$$
as is required. Now let $\cX$ be the $U/L$-wreath product. Then again by
the first part of~\eqref{121016t} for $\cY=\cX_{G/L}$, $M=G/L$, and $N=U/L$   we have
$$
\aut(\cX_{G/L})^{U/L}=\aut((\cX_{G/L})_{U/L})=\aut(\cX_{U/L}).
$$
Similarly, by the second part of~\eqref{121016t} for $\cY=\cX_U$, $M=U$, and $N=L$, we have 
$$
\aut(\cX_U)^{U/L}=\aut((\cX_U)_{U/L})=\aut(\cX_{U/L}).
$$
The two above equalities show that the hypothesis of Theorem~\ref{121016O} is satisfied and hence $\aut(\cX)^U=\aut(\cX_U)$. Applying the first of relations~\eqref{121016t} again for $\cY=\cX_U$, $M=U$, and $N=H$,
we get
$$
\aut(\cX)^H=(\aut(\cX)^U)^H=(\aut(\cX_U))^H=\aut(\cX_H),
$$
as is required.\eprf

\section{Extensions of coherent configurations and WL-algorithm}

\sbsnt{Partial order.}\label{300116a}
There is a natural partial order\, $\le$\, on the set of all coherent configurations on the same set~$\Omega$.
Namely, given two coherent configurations $\cX=(\Omega,S)$ and
$\cX'=(\Omega,S')$, we set
$$
\cX\le\cX'\ \Leftrightarrow\ S^\cup\subseteq (S')^\cup.
$$
The minimal and maximal elements with respect to this order are, respectively, the { trivial} and { complete} coherent configurations. 

\sbsnt{WL-algorithm.}
One of the most important properties of the partial ordering of coherent configurations comes from the fact that given a set $T\subseteq 2^{\Omega\times\Omega}$, there exists a unique minimal coherent configuration $\cX=(\Omega,S)$, for which $T\subseteq S^\cup$ (in particular, the set of all coherent configurations on the same set form a join-semilattice). This coherent configuration is called the {\it coherent closure} of $T$ and can be constructed by the well-known Weisfeiler-Leman algorithm (WL-algorithm) \cite[Section~B]{W76} in time polynomial in sizes of~$T$ and~$\Omega$. To stress this fact  the coherent closure of $T$
is denoted by $\WL(T)$.\medskip 

The {\it extension} of a coherent configuration $\cX$ with respect to the set $T$ is defined to be the coherent closure of $S\,\cup\,T$ and is denoted by~$\WL(\cX,T)$. The following statement is a straightforward consequence of~\cite[Theorem~8.2]{W76}.

\thrml{141016a}
Let $\cX=(\Omega,S)$ be a coherent configuration, $T\subseteq 2^{\Omega\times\Omega}$, and $\cY=\WL(\cX,T)$. Then 
$$
\aut(\cY)=\{f\in\aut(\cX):\ s^f=s\ \,\text{for all}\ s\in T\}.
$$ 
\ethrm

From Theorem~\ref{141016a}, it follows that if $X$ is a graph with vertex set~$\Omega$ and arc set~$R$, then $\aut(X)$ equals the automorphism group of the coherent closure $\WL(\{R\})$.

\sbsnt{Examples of extensions.}\label{131216a}
For a coherent configuration $\cX$ and an equivalence relation  $E\in\cE$, we denote by $\cX_E$ the extension of $\cX$ with respect to the set 
$$
T=\{1_\Delta:\ \Delta\in\Omega/E\}.
$$
From Theorem~\ref{141016a}, it follows that $\aut(\cX_E)$ consists of the automorphisms of $\cX$ leaving each class of~$E$ fixed. On the other hand, each automorphism of $\cX$ permute the classes of~$E$. Thus, $\aut(\cX_E)$ equals the kernel of the natural epimorphism from $\aut(\cX)$ to $\aut(\cX)^{\Omega/E}$.\medskip

In the above notation, let $c\in\sym(\Omega/E)$. Denote by $\cX_c$ the extension of $\cX$ with respect to the singleton $\{s\}$, where
\qtnl{150916u}
s=s(E,c)=\bigcup_{\Delta\in\Omega/E}\Delta\times\Delta^c
\eqtn
(note that if  $E=1_\Omega$, then the relation $s$ coincides with the graph of the permutation~$c$). Now by Theorem~\ref{141016a}, the group $\aut(\cX_c)$ consists of all $f\in\aut(\cX)$ such that $f^{\Omega/E}$ commute with~$c$.

\sbsnt{Extension of Cayley schemes.}
Let $\cX$ be a Cayley scheme over the group~$G$. Fix a group $H\in\fS$  and  consider the coherent configuration $\cY=\cX_E$, where $E=E_H$. Then the set $\Phi=\Phi(\cY)$ consists of the orbits of  $H$ acting on~$G$ by right multiplications. Besides, the equivalence relation $E$ is invariant with respect to the group $G_{right}$. Therefore, this group acts as an isomorphism group of the coherent configuration $\cY$. This enables us to define the coherent configuration
\qtnl{131216b}
\cY^G=(G,\{s^G:\ s\in S_\cY\}),
\eqtn
where $S_\cY$ is the set of basis relations of $\cY$ and $s^G$ is the union of the relations $s^g$,
$g\in G_{right}$. Obviously, 
$$
\cY\geq \cY^G\geq\cX.
$$ 
To find the group $\aut(\cY^G)$, we note that if $\rho:G_{right}\to \sym(\Phi)$ is the induced homomorphism, then 
$$
\im(\rho)=(G/H)_{right}\qaq\ker(\rho)\le\aut(\cY).
$$
In this situation, we can apply a result in \cite[Theorem~2.2]{E03} saying that 
\qtnl{121016k}
\aut(\cY^G)=G\aut(\cY).
\eqtn 

\lmml{121016i}
In the above notation, suppose that the group $G$ is abelian. Then $\cZ=\cY^G$ is a Cayley scheme over $G$. Moreover, 
\nmrt
\tm{1} $\aut(\cY)^H=\aut(\cZ)^H$ and $\cY_H=\cZ_H$,
\tm{2} $\cZ_{G/H}$ is the trivial cyclotomic scheme over $G/H$.
\enmrt
\elmm
\proof The fact that $\cZ$ is a Cayley scheme over~$G$ immediately follows from formula~\eqref{121016k}. This also implies statement~(2), because the coherent configuration $\cY^{\Omega/E}$ is complete. The aforementioned formula shows that the setwise stabilizer of the set $H$ in the group  $\aut(\cZ)$ coincides with $\aut(\cY)$. This proves the first equality of statement~(1). The second one follows from statement~(2) and \cite[Theorem~2.1]{EP10}.\eprf

\section{Structure of Cayley schemes over $E_4\times C_p$}\label{14102016r}

\sbsnt{Principal equivalence relation.} In this  section we are interested in the equivalence relations $E$ belonging to set $\cE=\cE(\cX)$, where $\cX=(\Omega,S)$ is a homogeneous coherent configuration on $n=4p$ points. The homogeneity of~$\cX$ implies that the valency $n_E$ of~$E$ divides~$n$. In what follows, we assume that $p\ge 5$ is a prime.

\dfntnl{120216t}
We say that $E\in\cE$ is a {\it principal} equivalence relation of~$\cX$ if one of the following statements hold:
\nmrt
\tm{E1} $n_E\geq p$ and $E$ is minimal (with respect to inclusion), 
\tm{E2} $n_E\leq 4$ and $E$  has a complement with respect to tensor or generalized wreath product.
\enmrt
In case {\rm (E2)}, it is assumed that $\cE$ contains no $E$ satisfying~{\rm (E1)}.
\edfntn

It follows from Definition~\ref{120216t} that every principal equivalence relation of the coherent configuration $\cX$ equals the equivalence closure of the union of at most two basis relations of~$\cX$. Since $|S|\le n$, we immediately obtain the following statement.

\lmml{240916b}
Given a homogeneous coherent configuration $\cX$ of degree $n=4p$, one can test in time $\poly(n)$ whether there exists a principal equivalence relation of~$\cX$, and (if so) find it within the same time.
\elmm

Clearly, for the trivial scheme of degree $n$, the equivalence relation $\Omega\times\Omega$ is principal. Let now $\cX$ be a Cayley scheme over a group~$G$ (of order~$n$), and let $E$ be a principal equivalence relation of~$\cX$.  If $\cX$ is a cyclotomic scheme, then the Sylow $p$-subgroup $P$ of $G$ is characteristic in $G$, and hence $P\in\fS$. Since also $|P|=p$, the equivalence relation $E_P$ satisfies the condition~(E1) and hence is a unique principal equivalence of~$\cX$.

\lmml{240916a}
Every Cayley scheme over the group~$G$ has a principal equivalence relation.
\elmm
\proof Let $\cX$ be a Cayley scheme over~$G$. By the above remarks, we may assume that the scheme $\cX$ is neither cyclotomic nor trivial. Therefore by Theorem~\ref{030916c}, the set $\fS$ contains two proper subgroups $U$ and $L$ of the group $G$ such that
\qtnl{121016g}
\cX=\cX_U\otimes\cX_L\qoq 
\cX\ \text{is the $U/L$-wreath product}.
\eqtn
Suppose that $\fS$ does not contain a principal equivalence relation satisfying condition~(E1). Then every minimal subgroup of $\fS$ is of order at most~$4$. It follows that if $|L|\le 4$, then $E_L$ is a principal equivalence satisfying condition~(E2), and we are done. Thus, we may assume that 
$$
L\ \,\text{is not minimal and}\ |L|>4.
$$ 
Then it is easily seen that $|L|=2p$. Now if the first equality in formula~\eqref{121016g} holds, then $E_U$ is a principal equivalence satisfying condition~(E2). In the remaining case, $L$ contains a proper subgroup, say~$H$, of order other than $p$. Therefore, $|H|=2$. Besides, $\cX$ being the $U/L$-wreath product is also the $U/H$-wreath product. Thus, $E_H$ is a principal equivalence satisfying condition~(E2).\eprf

\sbsnt{A principal equivalence relation of large valency.} 
Let $\cX$ be a Cayley scheme over the group $G$, $E$ a principal equivalence relation of $\cX$, and let $\cY=\cX_E$ be the coherent configuration defined in Subsection~\ref{131216a}.

\lmml{121016i1}
In the above notation, assume that $n_E\ge p$. Then for each set $\Delta\in G/E$, either $\aut(\cY)^\Delta=\sym(\Delta)$ or  $\cY_\Delta$ is isomorphic to a proper cyclotomic scheme over $\mZ_p$.
\elmm
\proof Let us verify that $E$ is minimal in $\cE(\cZ)$, where the Cayley scheme $\cZ:=\cY^G$ is defined by formula~\eqref{131216b}. First we note that 
$$
E\in\cE(\cX)\subseteq\cE(\cZ),
$$ 
because $\cZ$ is larger than~$\cX$. The minimality of $E$ in $\cE(\cZ)$ is obvious if $n_E=p$. In the other two cases, we have $n_E=2p$ or $4p$. This implies that $\cZ=\cX$ and the claim immediately follows from the definition of principal equivalence.\medskip 

Set $H=H_E$. We claim that to prove the lemma, it suffices to verify the validity of the equality
\qtnl{191016a}
\aut(\cZ)^H=\aut(\cZ_H).
\eqtn
Indeed, by Lemma~\ref{121016i} we have $\aut(\cY)^H=\aut(\cZ)^H$ and
$\cY_H=\cZ_H$. This implies that assuming \eqref{191016a}, 
$$
\aut(\cY)^H=\aut(\cZ)^H=\aut(\cZ_H)=\aut(\cY_H).
$$ 
By Corollary~\ref{12102016d} applied to~$\cX=\cZ$ and $G=H$, this implies that  the scheme $\cY_H$ is either trivial and then $\aut(\cY_H)=\sym(H)$, or is proper cyclotomic. This proves our claim, because the scheme $\cY_\Delta$ with $\Delta\in G/E$, is isomorphic to the scheme $\cY_H$ (the isomorphism is induced by any permutation of the group $G_{right}$ that takes $H$ to $\Delta$).\medskip

Let us prove equality~\eqref{191016a}. If the scheme $\cX$ is trivial, then we have $\cY=\cX$, $H=G$, and the equality is obvious. So we may assume that 
\qtnl{141216a}
|H|\in\{p,2p\}.
\eqtn
Then the scheme $\cZ\ge\cX$ is not trivial. If it is cyclotomic, then $|H|=p$ and the set $\fS(\cZ)$ contains the group $P\cong E_4$. Moreover, from statement~(2) of Lemma~\ref{121016i} it follows that $\cZ_P\cong\cZ_{G/H}$ is the trivial cyclotomic scheme over $G/H\cong P$. According to \cite[Lemma~2.3]{EKP13} this implies that $\cZ\cong\cZ_H\otimes\cZ_P$. Thus, by Theorem~\ref{030916c} the scheme $\cZ$ is proper tensor or generalized wreath product. Let us consider these two cases separately.\medskip

Let $\cZ=\cZ_U\otimes\cZ_L$, where $U$ and $L$ are proper subgroups of $G$ that belong to $\fS(\cZ)$. Without loss of generality, we may assume that $L$ is of order $p$ or $2p$. Then in view of~\eqref{141216a},
$$
(|H|,|L|)\in\{(p,p),(2p,p),(p,2p),(2p,2p)\}.
$$ 
By the minimality of $H$ the case $(2p,p)$ is impossible. By the same reason, in the case $(2p,2p)$ we have $H=L$, for otherwise $H\cap L$ is the proper subgroup of $H$. Thus, in any case $H\le L$. Therefore equality~\eqref{191016a} immediately follows from Lemma~\ref{12102016e}.\medskip

Let now $\cZ$ be the $U/L$-wreath product, where $U,L\in\fS(\cZ)$ are such that $1<L\le U<G$. By Lemma~\ref{12102016e}, to prove equality~\eqref{191016a} it suffices to verify that
\qtnl{191016q}
H\le L.
\eqtn
To this end, suppose first that $H\not\le U$. Then $E=E_H\not\subseteq E_U$ and hence the scheme~$\cZ$ has a basis relation $s\subseteq E\setminus E_U$. By the definition of generalized wreath product,  we have $E_L\cdot s=s$. It follows that $E_L$ is contained in the minimal equivalence relation $F\in\cE(\cZ)$ containing~$s$. Therefore $L\le H_F$. On the other hand, $H_F\le H$, because $s\subset E$. Thus, $L\le H_F\le H$. By the minimality of~$H$, this implies that $L=H$ which proves inclusion~\eqref{191016q}.\medskip

It remains to show that the statement
\qtnl{141216c}
H\le U\qaq H\not\le L
\eqtn
does not hold. Indeed, otherwise the minimality of $H$ implies that $H\cap L=1$ and hence $|L|=2$. Consequently, $\cZ'=\cZ_{G/H}$ is a Cayley scheme over the group $G'=G/H$ isomorphic to~$E_4$. This scheme is the $U'/L'$-wreath product with $U'=UH/H$ and $L'=LH/H$, because $\cZ$ is the $U/L$-wreath product. In view of assumption~\eqref{141216c}, 
$$
1<L'\le U'<G'.
$$
Therefore the rank of $\cZ'$ is less than~$4$. However, by  statement~(2) of Lemma~\ref{121016i}, the rank of $\cZ'$ equals~$4$, a contradiction.\eprf

\sbsnt{Summary}
The following theorem summarizes what we proved in this section and shows a way how we are going to use the principal equivalence relation.

\thrml{090916c}
Let $\cX=(\Omega,S)$ be a Cayley scheme over $E_4\times C_p$ ($p\ge 5$), $E$ a principal equivalence relation of~$\cX$, and $\cY=\cX_E$. Then
\nmrt
\tm{1} if $n_E\le 4$, then $\aut(\cX)^{\Omega/E}=\aut(\cX_{\Omega/E})$, 
\tm{2} if $n_E\ge p$, then one of the following statements holds:
\nmrt
\tm{a} $\cY$ is quasitrivial,
\tm{b} for each $\Delta\in\Phi(\cY)$, $\cY_\Delta$ is isomorphic to a proper cyclotomic scheme over~$\mZ_p$.
\enmrt
\enmrt
\ethrm
\proof If $n_E\le 4$, then the required statement immediately follows from the definition of the principal equivalence relation and Lemma~\ref{12102016e}. Otherwise, we are done by Lemma~\ref{121016i1} and the definition of a quasitrivial coherent configuration.\eprf

\section{Finding a representative set of semiregular $\mZ_p$-subgroups}\label{14102016}

In this section, we describe an efficient algorithm finding a {\it representative}  set~$B_p$ of the automorphism group~$K$ of a coherent configuration of degree~$4p$ (Subsection~\ref{100916d}). By definition, $B_p$ consists of semiregular $\mZ_p$-subgroups of~$K$ such that every group in $\reg(K,\mZ_p)$ is $K$-conjugate to one of them. Then we estimate the running time of the algorithm and explain the implementation details (Subsection~\ref{100916e}). In fact, the algorithm correctly finds $B_p$ if the input coherent configuration is isomorphic to a Cayley scheme over $G=E_4\times C_p$; we prove this in Subsection~\ref{100916f}. 

\sbsnt{The main algorithm.}\label{100916d}

In the algorithm below, we make use of the algorithm from~\cite{EP04} that constructs a {\it cycle base} of a coherent configuration~$\cX$; the cycle base is defined to be a maximal set of pairwise non-conjugated full cycles of the group $K=\aut(\cX)$. We always assume that $p\ge 5$ is a prime.\medskip

\centerline{\bf Main Subroutine (MS).}\medskip

\noindent {\bf Input:} A homogeneous coherent configuration~$\cX$ on $4p$ points.

\noindent {\bf Output:} A set $B_p\subseteq\reg(K,\mZ_p)$, where $K=\aut(\cX)$, that is empty or representative.
\medskip

\noindent{\bf Step 1.} Find a principal equivalence relation $E$ of~$\cX$; if there is no such~$E$, then output $B_p:=\varnothing$.\medskip

\noindent{\bf Step 2.} If $n_E\le 4$, then 

\def\asps{\hspace{-5mm}}
\nmrt
\item[] {\asps\bf Step 2.1.} find a cycle base~$C$ of the coherent configuration~$\cX_{\Omega/E}$;
\item[] {\asps\bf Step 2.2.} for each $c\in C$, find successively the relation $s=s(E,c)$ defined in~\eqref{150916u} and the coherent configuration $\cX_c=\WL(\cX,s)$; 
\item[] {\asps\bf Step 2.3.} output $B_p=\{P_c\in\reg(K,\mZ_p):\ c\in C\}$, where  $P_c$ is a Sylow $p$-subgroup of the group $K_c=\aut(\cX_c)$.\medskip
\enmrt

\noindent{\bf Step 3.} (Here $n_E\geq p$.) If the coherent configuration $\cY=\cX_E$ 
has a fiber $\Delta$ such that $|\Delta|\ne n_E$ or the homogeneous component $\cY_\Delta$ is not circulant, then output $B_p=\varnothing$.\medskip 

\noindent{\bf Step 4.} If the coherent configuration $\cY$ is quasitrivial, then output  $B_p=\{P\}$ with arbitrary $P\in\reg(\aut(\cY),\mZ_p)$.\medskip

\noindent{\bf Step 5.} Output $B_p=\reg(K',\mZ_p)$, where $K'\le\sym(\Omega)$ is the direct product of arbitrarily chosen regular cyclic subgroups $P_\Delta\le\aut(\cY_\Delta)$, $\Delta\in\Phi(\cY)$.

\sbsnt{Analysis of the running time.}\label{100916e} In the proof of the theorem below, we present detailed explanations of the steps of the algorithm MS and analyze the running time of each of them. 

\thrml{110916a}
The algorithm MS terminates  in time $\poly(n)$. In particular,
the size of its output is polynomially bounded.
\ethrm
\proof Let us verify successively that each step of the algorithm MS runs in polynomial time. \medskip

{\bf Step 1.} Here the required statement follows from Lemma~\ref{240916b}.\medskip

{\bf Step 2.} At Step 2.1, we apply the Main Algorithm from~\cite{EP04} that finds a cycle base of a coherent configuration in polynomial time. In our case, the input for this algorithm is the coherent configuration $\cX_{\Omega/E}$ of degree $p$ or $2p$. By the well-known upper bound for the size of a cyclic base (see~\cite{Mu99}), we have
\qtnl{110916g}
|C|\le |\Omega/E|\le 2p.
\eqtn

At Step 2.2, the coherent configuration $\cX_c$ found by the Weisfeiler-Leman algorithm is isomorphic to a coherent configuration extending the wreath product $\cY=\cX_1\wr\cX_2$, where $\cX_1$ is a trivial coherent configuration of degree $n_E=2$ or $4$, and $\cX_2$ is the trivial cyclotomic scheme over the group $\grp{c}$ of order $2p$ or $p$, respectively. 
Therefore,  the group~$K_c$ defined at Step~2.3 is isomorphic to a subgroup of the group 
\qtnl{141216r}
\aut(\cY)=\aut(\cX_1)\wr\aut(\cX_2)\cong\sym(n_E)\wr\mZ_{n/n_E},
\eqtn
which is solvable and can be constructed efficiently. This enables to construct in polynomial time the group $K_c$ by the Babai-Luks algorithm~\cite{BL} and then its Sylow $p$-subgroup $P_c$ by the Kantor algorithm (see~\cite{S03}). Thus, in view of inequality~\eqref{110916g}, Step~2 terminates in time $\poly(n)$.\medskip

{\bf Step 3.} This step involves the WL-algorithm to construct
the coherent configuration~$\cY$ and the Main Algorithm from~\cite{EP04} to
find a cycle base $C_\Delta$ of the group $\aut(\cY_\Delta)$: the coherent configuration $\cY_\Delta$ is circulant if and only if the set $C_\Delta$ is
not empty. Since both of algorithms run in time $\poly(n)$, we are done.\medskip 

{\bf Step 4.} To verify whether the coherent configuration~$\cY$ is quasitrivial, it suffices to make use of Lemma~\ref{090916a}. Next, assume that $\cY$ is quasitrivial. Then  one can efficiently find the decomposition~\eqref{091016a} of the coherent configuration $\cY$ into the direct sum  and enumerate the fibers of $\cY$ so that
$$
\Phi(\cY)=\{\Delta_{ij}:\ i=1,\ldots,m,\ j=1,\ldots,m_i\}.
$$
In our case, $1\le m\le 4$ and $\sum_im_i=n/n_E$. For each $i$, one can efficiently find  a permutation  $g_i\in\sym(\Delta_{i1})$, which is the product of $n_E/p$ disjoint $p$-cycles. Then applying the second part of Theorem~\ref{090916b} for $\cX=\cY$, we see that the permutation $g$ defined by
formula~\eqref{101016a} is an automorphism of the coherent configuration~$\cY$. In particular,
the group $P=\grp{g}$ is contained in $\reg(\aut(\cY),\mZ_p)$.\medskip

{\bf Step 5.} Here, we define $P_\Delta$ as the group generated by an arbitrary full cycle belonging to the set $C_\Delta$ found at Step~3. Then 
$$
|K'|=\prod_{\Delta\in\Omega/E}|P_\Delta|=|\Delta|^{|\Omega/E|}\le n^4.
$$
Thus, all the elements of order~$p$  in $K'$ and hence the set $B_p$ can be found in time $\poly(n)$.\eprf

\sbsnt{The correctness of the MS}\label{100916f} We keep the notation of Subsection~\ref{100916d} and denote by~$\cK_p$ the class of coherent configurations isomorphic to a Cayley scheme over a group $G=E_4\times\mZ_p$. Set
$$
\reg_p(K)=\{H_p:\ H\in\reg(K,G)\},
$$
where $H_p$ is the Sylow $p$-subgroup of the group~$H$.

\thrml{030916d}
Let $\cX$ be a coherent configuration on $4p$ points, and let the group $K$ and set $B_p$ be as in the Main Subroutine. Then if $\cX\in\cK_p$, then $B_p$ is not empty and every group in $\reg_p(K)$ has a $K$-conjugate in $B_p$. In particular, the set $B_p$ is representative.
\ethrm

\crllrl{150916t}
If $B_p=\varnothing$, then $\cX\not\in\cK_p$.
\ecrllr

{\bf Proof of Theorem~\ref{030916d}.} Without loss of generality, we may assume that $\cX$ is a Cayley scheme over $G=E_4\times\mZ_p$. Then by Lemma~\ref{240916a}, the set $\cE=\cE(X)$ contains a principal equivalence relation~$E$. Therefore, at Step~1 the algorithm MS does not terminate. 
Let us verify that every group 
$$
Q\in\reg_p(K)
$$ 
has a $K$-conjugate in $B_p$. To this end, let $H\in \reg(K,G)$ be such that $Q=H_p$.\medskip

Assume first that $n_E\le 4$. Then $n_E=2$ or $4$. Therefore, $H^{\Omega/E}$ is a regular cyclic group of order~$n/n_E=rp$ with $r=1$ or $2$. Consequently, there exists a permutation $h\in H$ of order $rp$ such that  $H^{\Omega/E}$ is  generated by $h^{\Omega/E}$. In particular, 
$$
Q=H_p=\grp{h}_p.
$$ 
On the other hand, $K^{\Omega/E}=\aut(\cX_{\Omega/E})$ by statement~(1) of Theorem~\ref{090916c}. Now, if $C$ is the cycle base found at Step~2.1, then there exist a  permutation $k\in K$ and a full cycle $c\in C$ such that
$$
(h^k)^{\Omega/E}=(h^{\Omega/E})^{k^{\Omega/E}}=c.
$$
By the definition of the relation $s$ at Step~2.2, this immediately implies that $s^{h'}=s$, where $h'=h^k$. Since also $h'\in K$, we conclude that $h'$ belongs to the group~$K_c$  constructed at Step 2.3. However, the order of $h'$ equals $rp$ (the order of $h$) and the order of $P_c$ is equal to~$p$ (see formula~\eqref{141216r}). Thus, there exists $k'\in K_c$
such that $(h'_p)^{k'}\in P_c$, where $h'_p$ is the $r$th power of~$h'$. 
Now,
$$
(Q)^{kk'}=(\grp{h}_p)^{kk'}=\grp{(h')^{k'}}_p=\grp{(h'_p)^{k'}}=P_c.
$$
Thus, the set $B_p$ constructed at Step 2.3 contains the $K$-conjugate $P_c$ of the group $Q$, as required.\medskip

Assume that $n_E\ge p$. Then $n/n_E$ is less or equal than~$4$. Since $p\ge 5$, we conclude that $Q$ acts trivially on $\Omega/E$. It follows that 
\qtnl{141216y}
Q\le\aut(\cY),
\eqtn
where $\cY$  is the coherent configuration found at Step~$3$. Since $\cX$ is a Cayley scheme over~$G$, every  fiber $\Delta\in\Phi(\cY)$ is a class of the equivalence relation~$E$, and hence $|\Delta|=n_E$. The group $\aut(\cY_\Delta)$ contains a regular subgroup $G^\Delta$. This group is cyclic if $n_E\ne 4p$ and hence $\cY_\Delta$ is a circulant scheme in this case.
Finally, $n_E=4p$, then $\cY_\Delta$ is trivial and hence circulant. Thus, the algorithm MS does not terminate at Step~3.\medskip

Suppose first that the coherent configuration $\cY$ is quasitrivial. Then in the notation of Theorem~\ref{090916b}, we have
$$
Q^{\Delta_{i1}}\le\aut(\cY)^{\Delta_{i1}}=\sym(\Delta_{i1}),\quad i=1,\ldots,m.
$$
However,  $Q^{\Delta_{i1}}$ and $P^{\Delta_{i1}}$ are semiregular cyclic subgroups of order $p$ in the group $\sym(\Delta_{i1})$, where $P$ is the group defined at Step~4. Therefore, there exists a permutation $g_i\in\sym(\Delta_{i1})$ such that 
$$
(Q^{\Delta_{i1}})^{g_i}=P^{\Delta_{i1}},\quad i=1,\ldots,m.
$$
Then $Q^g=P$, where the permutation $g\in K$ is defined by formula~\eqref{101016a}. By Theorem~\ref{090916b}, we have $g\in \aut(\cY)$. Since the latter group is contained in~$K$,
the set $B_p$ constructed at Step 4 consists of the $K$-conjugate $P$ of the group $Q$, as required.\medskip
 
To complete the proof, we may assume that $\cY$ is not quasitrivial. Then
by statement~(2) of Theorem~\ref{090916c},  for each $\Delta\in\Phi(\cY)$ the coherent configuration $\cY_\Delta$ is isomorphic to a proper cyclotomic scheme over~$\mZ_p$. In particular, we come to Step~5 with $|\Delta|=p$. Moreover, 
$$
\reg(\aut(\cY_\Delta),\mZ_p)=\{P_\Delta\},
$$
where $P_\Delta$ is the group found at Step~5. By iclusion~\eqref{141216y}, this implies  that 
$$
\orb(Q,\Omega)=\Phi(\cY).
$$
Therefore, $Q^\Delta\in\reg(\aut(\cY_\Delta),\mZ_p)$  and hence $Q^\Delta=P_\Delta$ for all $\Delta$. This implies that 
$$
Q\le \prod_{\Delta\in\Omega/E} Q^\Delta= \prod_{\Delta\in\Omega/E}P_\Delta=K',
$$
where $K'$ is the group defined at Step~5. Thus, $Q$ is contained in the set~$B_p$ found at this step and we are done.\eprf

\section{Proof of Theorem~\ref{030916a}}\label{100916g}
Let $\cX$ be a coherent configuration constructed from the graph $X$ by the WL-algorithm. Then $X$ is isomorphic to a Cayley graph over $G$ if and only if $\cX$ is isomorphic to a Cayley scheme over~$G$. Since also $\aut(X)=\aut(\cX):=K$, it suffices to check in time $\poly(n)$ whether the set $\reg(K,G)$ is not empty, and (if so) find an element of this set within the same time.\medskip

Without loss of generality, we assume $p\ge 5$. Let $B_p\subseteq\reg(K,\mZ_p)$ be the set constructed by the algorithm MS applied to the coherent configuration~$\cX$. By Theorem~\ref{110916a} this can be done in time $\poly(n)$. If the set $B_p$  is empty, then $\cX\not\in\cK_p$ by Corollary~\ref{150916t} and hence the set $\reg(K,G)$ is also empty. Thus, we may assume that $\cX\in \cK_p$ and hence $B_p\ne\varnothing$. For each group $P\in B_p$, we define a set  
$$
R(P)=\{H\in\reg(C_P\cap K,G):\ P\le H\}
$$
of regular subgroups of $\sym(\Omega)$, where $C_P$ is the centralizer of $P$ in $\sym(\Omega)$. Note that the group $C_P$ is permutation isomorphic to the wreath product $\mZ_p\wr\sym(4)$. Therefore, 
\qtnl{101016w}
|C_P|=|\mZ_p\wr\sym(4)|=24p^4\le n^4,
\eqtn
and the group $C_P\cap K$ can be found in time $\poly(n)$ by testing each permutation of $C_P$ for membership to the group~$K$. Next, every group $H\in R(P)$ is generated by $P$ and two involutions $x,y\in C_P\cap K$. Since the number of such pairs $(x,y)$ does not exceed $|C_P|^2$, the following statement is a consequence of inequality~\eqref{101016w}.

\lmml{100916x}
Given a group $P\in B_p$, the set $R(P)$ can be found in time $\poly(n)$. In particular, $|R(P)|\le n^c$ for a constant $c>0$.
\elmm

This lemma shows that to complete the proof we need to verify the implication
$$
\reg(K,G)\ne\varnothing\quad\Rightarrow\quad \bigcup_{P\in B_p}R(P)\ne\varnothing.
$$ 
To this end, it suffices to prove that every  group $V\in\reg(K,G)$ is $K$-conjugate to a group belonging to $R(P)$ for some $P\in B_p$. However by the above assumption,  $\cX\in\cK_p$ and $V_p\in\reg_p(K)$, where $V_p$ denotes the Sylow $p$-group of $V$. Therefore, by Theorem~\ref{030916d}, there exists $P\in B_p$ and $k\in K$ such that $(V_p)^k=P$. It follows that 
$$
V^k\le (C_K(V_p))^k=C_K((V_p)^k)=C_K(P)\le C_P,
$$
where $C_K(V_p)$ denotes the centralizer of $V_p$ in $K$.
Since also $V^k\le K$ and $V^k\cong G$, we conclude that $V^k\in R(P)$,
as required. 

\section*{Acknowledgement}
The paper was finished during the research stay of the second author at the Faculty of Applied Sciences of the University of West Bohemia in October, 2016 supported by the project P202/12/G061 of Czech Science Foundation. The first author is supported by the grants APVV-15-0220, VEGA 1/0150/14, Project LO1506 of the Czech Ministry of Education, Youth and Sports and Project P202/12/G061 of Czech Science Foundation.


\begin{thebibliography}{99}
	
\bibitem{B77}
L.~Babai, {\em Isomorphism problem  for a class of point symmetric structures}, Acta Math. Acad. Sci. Hung., {\bf 29}, No.~3, 329--336 (1977).
	
\bibitem{BL}
L.~Babai and E.~M.~Luks, {\em Canonical labeling of graphs}, Proc. 15th ACM
STOC, 171--183 (1983).
	

\bibitem{E03} 
S.~Evdokimov and I.~Ponomarenko, {\em Characterization of cyclotomic schemes and normal Schur rings over a cyclic group},  St.~Petersburg Math. J.,  {\bf 14},  No.~2, 189--221  (2003).

\bibitem{EP04} 
S.~Evdokimov and I.~Ponomarenko, {\em Circulant graphs: recognizing and isomorphism testing in polynomial time}, St. Petersburg Math. J., {\bf 15}, 813--835 (2004).

\bibitem{EP05}
S.~Evdokimov and I.~Ponomarenko, {\em A new look at the Burnside-Schur 	theorem}, Bulletin of the London Mathematical Society, {\bf 37} (2005), 535--546.

\bibitem{EP09}
S.~Evdokimov and I.~Ponomarenko, {\em Permutation group approach to association schemes}, European J. Combin., {\bf 30}, No. 6, 1456--1476 (2009).

\bibitem{EP10}
S.~Evdokimov and I.~Ponomarenko, {\em Schemes of a finite projective plane and their extensions}, St. Petersburg Math. J., {\bf 21}, No.~1, 65--93 (2010).

\bibitem{EP12}
S.~Evdokimov and I.~Ponomarenko, {\em Schurity of S-rings over a cyclic group and generalized wreath product of permutation groups}, St.
Petersburg Math. J., {\bf 24}, No.~3, 431--460 (2013).


\bibitem{EKP13} 
S.~Evdokimov, I.~Kov\'acs, and I.~Ponomarenko, {\em On schurity of finite abelian groups}, Commun. Algebra, {\bf 44}, No.~1, 101--117 (2016).

\bibitem{KP}  
M.~Klin and R.~P\"oschel, {\em The K\"onig problem, the isomorphism problem for cyclic graphs and the method of Schur rings}, In: ``Algebraic Methods in Graph Theory, Szeged, 1978'', Colloq. Math. Soc. J\'anos Bolyai, Vol. 25, North-Holland, Amsterdam (1981), pp. 405--434.

\bibitem{KM09}
I.~Kov\'acs and  M.~Muzychuk, {\em The group $\mZ_p^2\times\mZ_q^{}$ is a CI-group}, Commun. Algebra, {\bf 37}, No.~10, 3500--3515 (2009).

\bibitem{L16}
A.~Lang, {\em A Classification of the Supercharacter Theories of $C_p\times C_2\times C_2$ for prime $p$}, {\tt arXiv:1609.07182[math.RT]}, 1--28 (2016).

\bibitem{Mu99}
M.~Muzychuk, {\em On the isomorphism problem for cyclic combinatorial 	objects}, Discr. Math., {\bf 197/198}, 589--606 (1999).

\bibitem{MP} 
M.~Muzychuk and I.~Ponomarenko, {\em Schur rings}, European J. Combin., {\bf 30}, No.~6, 1526--1539 (2009).

\bibitem{S03}
A.~Seress, {\em Permutation group algorithms}, Cambridge Tracts Math., 
{\bf 152}, Cambridge University Press, Cambridge (2003).

\bibitem{W76}
B.~Weisfeiler (editor), {\em On construction and identification of graphs}, Lecture Notes Math., {\bf 558} (1976).

\bibitem{W64}
H.~Wielandt, {\em Finite permutation groups}, Academic Press, New York - London (1964).


\end{thebibliography}
\end{document}